\documentclass[ 12pt]{article}
\usepackage[utf8]{inputenc}
\usepackage[left=1.2in,
            right=1in,
            top=1in,
            bottom=1in]{geometry}

\usepackage{amsfonts}
\usepackage{hyperref}
\usepackage{amssymb}
\usepackage{amsmath}
\usepackage{amsthm}

\usepackage{graphicx}

\usepackage{algorithm,algorithmic}

\newtheorem{theorem}{Theorem}

\newtheorem{lemma}[theorem]{Lemma}

\RequirePackage{epsfig}
\RequirePackage{amssymb}
\RequirePackage{natbib}
\RequirePackage{graphicx}

\if@usehyper
\RequirePackage[colorlinks=false,allbordercolors={1 1 1}]{hyperref}
\fi

\bibpunct{(}{)}{;}{a}{,}{,}

\def\jmlrheading#1#2#3#4#5#6{\def\ps@jmlrtps{\let\@mkboth\@gobbletwo%
\def\@oddhead{\scriptsize Journal of Machine Learning Research #1 (#2) #3 \hfill Submitted #4; Published #5}%
\def\@oddfoot{\scriptsize \copyright #2 #6. \hfill}%
\def\@evenhead{}\def\@evenfoot{}}%
\thispagestyle{jmlrtps}}

\def\ShortHeadings#1#2{\def\ps@jmlrps{\let\@mkboth\@gobbletwo%
\def\@oddhead{\hfill {\small\sc #1} \hfill}%
\def\@oddfoot{\hfill \small\rm \thepage \hfill}%
\def\@evenhead{\hfill {\small\sc #2} \hfill}%
\def\@evenfoot{\hfill \small\rm \thepage \hfill}}%
\pagestyle{jmlrps}}

\usepackage{natbib}
\hypersetup{
    colorlinks,
    breaklinks,
    linkcolor = blue,
    citecolor = blue,
    urlcolor  = blue,
}
\usepackage{url}            
\usepackage{graphicx,subfigure,color}
\usepackage{amsmath,mathrsfs,amssymb,hyperref}
\usepackage{xcolor}

\usepackage{hyperref}

\title{Improving Dynamic Regret in Distributed Online Mirror Descent Using Primal and Dual Information}
\author{\textbf{Nima Eshraghi} \hspace{4cm}\hfill {neshraghi@ece.utoronto.ca} 
   \\ \textbf{Ben Liang}       \hspace{4cm} \hfill    {liang@ece.utoronto.ca}\\
    Department of Electrical and Computer Engineering,  \\
   University of Toronto, Toronto, Canada}

\date{}

\begin{document}

\maketitle

\begin{abstract}
We consider the problem of distributed online optimization, with a group of learners connected via a dynamic communication graph. The goal of the learners is to track the global minimizer of a sum of time-varying loss functions in a distributed manner. We propose a novel algorithm, termed Distributed Online Mirror Descent with Multiple Averaging Decision and Gradient Consensus (DOMD-MADGC), which is based on mirror descent but incorporates multiple consensus averaging iterations over local gradients as well as local decisions. The key idea is to allow the local learners to collect a sufficient amount of global information, which enables them to more accurately approximation the time-varying global loss, so that they can closely track the dynamic global minimizer over time. We show that the dynamic regret of DOMD-MADGC is upper bounded by the path length, which is defined as the cumulative distance between successive minimizers. The resulting bound improves upon the bounds of existing distributed online algorithms and removes the explicit dependence on $T$.
\end{abstract}

\section{Introduction}
In recent years, applications have emerged that require extremely large data volumes over many networked machines. As a result, distributed collection and processing of these datasets are not only desirable but often necessary \citep{boyd2011distributed}.
Therefore, distributed optimization has become popular in solving problems that arise in various areas of control and learning \citep{boyd2011distributed, duchi12dual, xi5526distributed, rabbat15}.

The aforementioned distributed optimization methods assume a static loss function. Nevertheless, in many practical applications, the system parameters and loss functions vary over time. For example, time-varying loss functions frequently appear in online machine learning, where data samples arrive dynamically, so that newly observed data samples result in new losses.  Another example is object-tracking, where the goal is to track the time-varying states of a moving target. These problems can be solved using online optimization algorithms that update the decisions based on the dynamically arriving data.

The performance of online optimization algorithms is usually measured in terms of regret. Depending on the problem environment and settings, different notions of regret have been proposed in the literature. For instance, static regret measures the accumulated loss incurred by the online optimization algorithm against that of an offline optimal solution, which is made in hindsight given the knowledge of all loss functions.  The static regret of online optimization algorithms has been extensively studied in the literature \citep{cesa2006prediction, shalev2012online,  hazan2019introduction, orabona2019modern}. 

Dynamic regret is a more appropriate performance metric when the underlying target of interest is time-varying. Dynamic regret measures the accumulated loss experienced by the online algorithm versus that of a sequence of optimizers, which separately minimize the loss functions at every round. In contrast to the static regret, establishing a sublinear upper bound on the dynamic regret is generally impossible, due to the arbitrary and potentially severe fluctuations in the problem environment. Instead, the growth of the dynamic regret depends on some regularity measures, which reflect the speed of changes in a dynamic environment. For instance, the \textit{path length} is a popular regularity measure that is related to the rate of changes in the sequence of minimizers over time, i.e.,\vspace*{-0.4cm} 
\begin{align}\label{q_pl_def}
    C_T = \sum_{t =2}^{T} \| x^*_t - x^*_{t-1} \|, 
\end{align}
where $x^*_t = {\text{argmin}_{x \in \mathcal{X}}} f_t(x)$ is the minimizer of the time-varying loss function. When loss functions are convex, the dynamic regret of online convex optimization algorithms can be bounded by $O(\sqrt{T} (C_T+1))$ \citep{zinkevich2003online}, which can be improved to $O(\sqrt{T (C_T+1)})$ when knowledge of the path length $C_T$ and the number of rounds $T$ is available in advance \citep{yangb16}. 
Recent studies have shown that the dependency of dynamic regret on the number of online rounds can be removed when stronger assumptions on the function curvature is available \citep{mokhtari2016online, zhang2017improved}. In particular, a dynamic regret bound of $O(1 + C_T)$ has been derived for strongly convex and smooth loss functions \citep{mokhtari2016online, zhang2017improved}, which improves the prior bound of \citep{zinkevich2003online} by removing a $\sqrt{T}$ factor. 

However, \citep{mokhtari2016online, zhang2017improved} focus on centralized problems, where there is only a single learner. Distributed optimization can be more challenging since the learners have only partial information about the global problem, which necessitates engagement in communication so that they can complement their insufficient knowledge. None of the previously proposed distributed online convex optimization algorithms \citep{shahrampour2017distributed, dixit2019distributed, tarokh2020distributed} achieve $O(1+C_T)$ dynamic regret.

In this paper, we investigate whether it is possible to obtain the aforementioned dynamic regret bound of $O(1+ C_T)$ in distributed online convex optimization. To this end, we propose a new algorithm termed Distributed Online Mirror Descent with Multiple Averaging Decision and Gradient Consensus (DOMD-MADGC), which is a variant of DOMD with an improved dynamic regret. Previous works on distributed mirror descent \citep{shahrampour2017distributed, xi5526distributed, rabbat15} perform a single consensus step only on the local decision vectors, and then each learner tries to minimize its local loss using its local gradient information. The resulting lack of knowledge about the global gradient prevents the learners from closely tracking the time-varying global minimizer. In contrast, DOMD-MADGC allows the communication of the local gradient vectors in addition to the local decisions, leading to reduced global loss. 
We further observe that simply averaging the local decisions and the local gradients in a single consensus step does not provide the dynamic regret bound that we seek. Therefore, we design DOMD-MAGDC to use multiple consensus averaging iterations per online round to collect the local information update from distant learner nodes. Our analysis reveals that a logarithmically increasing number of consensus iterations will result in sufficiently fast decrease of local errors. With the proposed design, DOMD-MADGC achieves a dynamic regret bound of $O(1+C_T)$, without any prior knowledge of the path length $C_T$ or the number of round $T$.

\section{Related Works} \label{related_works}

Early works on online convex optimization mostly concern static regret. In the centralized setting, it has been established that various optimization methods including online gradient descent \citep{zinkevich2003online}, online dual averaging \citep{xiao2010dual}, online mirror descent \citep{duchi2010composite}, and many others \citep{ shalev2012online,  hazan2019introduction} achieve an upper bound of $O(\sqrt{T})$ and $O(\log T)$ on the static regret, for convex and strongly convex loss functions, respectively. The static regret of distributed online convex optimization algorithms have also been extensively studied in the literature \citep{hosseini2013online, mateos2014distributed, akbari2015distributed, tsianos2016efficient, lee2016coordinate, yuan2020distributed}, where the same regret rates have been derived under similar convexity assumptions.   However, it is not previously known whether similar results hold for the more useful dynamic regret. That is, whether there exists a distributed online algorithm that has dynamic regret bound comparable to the most competitive centralized online algorithms.

There has been a surge of recent interest in studying dynamic regret in the centralized setting \citep{hall2013dynamical, hall2015, jadbabaie2015online, yangb16, mokhtari2016online, zhang2017improved, chang2020unconstrained}. Previous works on online gradient descent \citep{zinkevich2003online}, and online mirror descent \citep{hall2013dynamical, hall2015} show that the dynamic regret of both methods is bounded by $O(\sqrt{T}(1 + C_T))$ when the loss functions are convex, which can be tightened to $O(\sqrt{T(1 + C_T)})$ when the knowledge of path length $C_T$ is present a priori \citep{yangb16}. Furthermore, when the loss functions are strongly convex and smooth, the upper bound on the dynamic regret can be improved to $O(1+C_T)$ \citep{mokhtari2016online, zhang2017improved}. In another direction, researchers have used various regularity measures to bound the dynamic regret of different online optimization algorithms \citep{jadbabaie2015online, besbes2015non, chang2020unconstrained, campolongo2020temporal, campolongo2021closer}. The above works do not apply to the distributed setting.

Distributed online mirror descent is studied in \citep{shahrampour2017distributed}, where a dynamic regret bound of $O(\sqrt{T(1+C_T)})$ is derived when the joint knowledge of the path length $C_T$ and the number of round $T$ is available in advance. However, it is generally impossible to know $C_T$ beforehand in practice, so most subsequent works do not require this assumption. A distributed online gradient tracking algorithm is proposed in \citep{lu2019online}, which has a dynamic regret bound of $O( \sqrt{1+C_T} T^{3/4} \sqrt{\ln T})$.  The dynamic regret of distributed online proximal gradient descent with an $O(\log t)$ communication steps per round is bounded by $O(\log T (1+ C_T))$ \citep{dixit2019distributed}. 
Using the gradient variation as the regularity measure $V_T = \sum_{t=2}^T \sup_{x \in \mathcal{X}}\|\nabla f_t(x) - \nabla f_{t-1}(x)\|$, a dynamic regret bound of $O(\sqrt{T}(1+C_T) + V_T + \sqrt{V_T C_T})$ is derived in \citep{li2020distributed}, which is improved to $O(1 + C_T + V'_T)$ in \citep{tarokh2020distributed}, where $V'_T$ is a variant of $V_T$ in which gradients are evaluated at the optimal points. However, the work of \citep{tarokh2020distributed} studies only unconstrained online optimization problem.

In this work,  we show that the proposed DOMD-MADGC has a dynamic regret bound of $O(1+C_T)$. In contrast to most prior works, where the dynamic regret bound depends on a combination of a sublinear term in $T$ (e.g., $\log T$ or $\sqrt{T}$) and some regularity measure, our bound depends only on the path length $C_T$.  Furthermore, the resulting bound of $O(1+C_T)$ improves upon the existing regret bounds in dynamic distributed learning. Unlike \citep{shahrampour2017distributed}, the dynamic regret bound of our proposed algorithm can be achieved by selecting a step size without the knowledge of the path length $C_T$ or the number of rounds $T$. Unlike \citep{tarokh2020distributed}, in this work, we allow the optimization variables to be constrained by an arbitrary compact and convex set.

\section{Problem Formulation}

\subsection{Network Model}

We consider a distributed learning network in a time-slotted setting. The learner nodes interact with each other over a time-varying  communication topology, which is modeled by an undirected graph $\mathcal{G}_t = (\mathcal{V}, \mathcal{E}^t)$, where  $\mathcal{V} = \{1, \ldots, n\}$ and  $\mathcal{E}^t \subset \mathcal{V} \times \mathcal{V}$ denote the set of nodes and edges present at round~$t$, respectively. 

Different from online convex optimization where only a single learner is present, in distributed online optimization, every node $i \in \mathcal{V}$ is a local learner that receives a private local loss in every round. Each learner node $i$ is associated with a sequence of time-varying loss functions, i.e., $\{f_{i,1}(x), f_{i,2}(x), \ldots, f_{i,T}(x)\}$, where $f_{i, t}$ denotes the loss function of node $i$ at round $t$, and  $x \in \mathcal{X}$ represents a decision variable taken from a compact and convex set $\mathcal{X}$.  At each round $t$, after each learner $i$ submits a  decision $x_{i,t}$, the local loss function $f_{i,t}(x)$ is revealed, and the learner suffers a corresponding loss of $f_{i, t}(x_{i, t})$. Each local loss is only observed by a single learner. Therefore, the learners need to interact with each other, to supplement their incomplete knowledge of the global task.

\subsection{Distributed Online Optimization}

We are interested in an optimization problem with a global loss function, represented by $f_t(x)$ at round $t$. It is based on the local loss functions that are distributed over the entire  network, i.e.,
\begin{align} \label{q_glob_loss}
f_t(x) = \sum_{ i = 1}^{n} f_{i, t}(x).
\end{align}

The goal of the online learners is to minimize the total loss over a finite number of rounds~$T$. Prior studies on online learning often measure the quality of decisions in terms of \textit{static} regret, defined as the difference between the total loss incurred by the online algorithm and that of an optimal fixed offline solution, which is made in hindsight with knowledge of $f_t(x)$ for all $t$. The benchmark variable in static regret is a fixed point that resides in the feasible set $\mathcal{X}$. Such a static metric can accurately reflect the performance of an online algorithm as long as the static benchmark performs consistently well over all rounds. However, this does not always hold in a dynamic environment, where the sequence of loss functions are time-varying. Thus, \textit{dynamic} regret has been proposed as a more stringent metric that incorporates a time-varying comparator sequence. Most commonly, the performance of online algorithms is measured relative to a sequence of minimizers \citep{shahrampour2017distributed, tarokh2020distributed}, i.e,
\begin{align} \label{def_dy_reg}
\mathrm{Reg}_T^d = \frac{1}{n} \sum_{i =1}^{n}\sum_{t=1}^{T}f_t(x_{i,t}) - \sum_{t=1}^{T}f_t(x^*_t).
\end{align}
where $x^*_t = {\text{argmin}_{x \in \mathcal{X}}} f_t(x)$ is a minimizer of the global loss function at round $t$.

It is well-known that the online optimization problems may be intractable in a dynamic setting, due to arbitrary fluctuation in the loss functions. Therefore, obtaining a sublinear bound on the dynamic regret may be impossible. However, it is possible to bound the dynamic regret in terms of regularity measures that reflect how fast an environment changes as time progresses.  A popular quantity to represent such regularity is the path length $C_T$, defined in \eqref{q_pl_def}, which collects the variation between successive minimizers.

Since the local loss function is available only to the local learner, due to partial access to the global information, it is unlikely for each learner to accurately compute the global minimizer $x^*_t$ at each round $t$. To minimize the dynamic regret, a careful design of the online algorithm is required to allow distributed learners to collect sufficient global information to track the time-varying global minimizer. Thus, it is necessary to use the local information as well as those collected from the neighboring learners, so that each learner can find its local estimates sufficiently close to the global gradient $\nabla f_t(x^*_t)$ and minimizer $x^*_t$. 

\subsection{Preliminaries}\label{pre}

To be self-contained, we collect here several standard definitions that are used in this paper.

\textit{Definition 1}: A function $f(x)$ is Lipschitz continuous with factor $G$ if for all $x$ and $y$ in $\mathcal{X}$, the following holds:
\begin{align} \nonumber
 | f(x) - f(y) | \leq G \| x - y\|, ~\forall x,y \in \mathcal{X}.
\end{align}

\textit{Definition 2}: A subdifferentiable function $f(x)$ is $\beta$-smooth with respect to some norm $\|\cdot\|$, if there exists a positive constant $\beta$ such that
\begin{align} \nonumber
f(y) \leq f(x) + \langle \nabla f(x), y - x \rangle + \frac{\beta}{2} \| y - x \|^2, ~\forall x,y \in \mathcal{X}.
\end{align}
where $\nabla f(x)$ stands for the subgradient of function $f(x)$.

\textit{Definition 3}: A subdifferentiable convex function $f(x)$ is $\lambda$-strongly convex with respect to some norm $\|\cdot\|$, if there exists a positive constant $\lambda$ such that
\begin{align} \nonumber
 f(y) + \langle \nabla f(y), x - y \rangle + \frac{\lambda}{2} \| x - y \|^2 \leq f(x), ~\forall x,y \in \mathcal{X}.
\end{align} 

\textit{Definition 4}: The Bregman divergence with respect to the regularization function $r(x)$ is defined as 
\begin{align} \nonumber
D_{r}(x,y) = r(x) - r(y) - \langle \nabla r(y), x - y \rangle.
\end{align}

The Bregman divergence is a general distance-measuring function, which contains the Euclidean norm and the Kullback-Leibler divergence as two special cases.

\section{Distributed Online Mirror Descent with Multiple Averaging Decision and Gradient Consensus}
\label{alg_desc}

We now present DOMD-MADGC, a distributed online optimization algorithm that incorporates multiple consensus averaging of local gradients as well local decisions, over time-varying dynamic learning networks. We will bound the performance of DOMD-MADGC in terms of  dynamic regret, showing that it achieves the same regret rate as the centralized algorithms of \citep{mokhtari2016online, zhang2017improved}.

\subsection{Algorithm Description}

DOMD-MADGC uses mirror descent as its optimization engine.   The core of the algorithm is the use of consensus averaging over both the local decisions and the local gradients. 
We note that in prior studies on mirror descent, only the local decisions are used to reach a consensus, while DOMD-MADGC applies consensus averaging to the local gradients as well. We motivate this design choice by comparing it with standard online mirror descent, where each learner performs consensus averaging over the local decision vectors, and then uses its own local gradient to update its decision. Therefore, in standard online mirror descent, while each learner aims to stay close to the  averaged decision, it takes a step in a direction to minimize its local loss, which may not align well with the direction toward the global minimizer. In other words, each local learner has access to only one component of the global loss, and the resulting lack of knowledge may prevent the learner from closely tracking the time-varying global target. Therefore, DOMD-MADGC uses both the local gradients and the local decisions in the consensus phase.
As a result, the learners can use a local approximation of the global gradient in their decision updates, which leads them to find decisions that reduces the global loss function. 

However, simply adding gradient consensus to standard online mirror descent still does not sufficiently reduce the dynamic regret. The success of the algorithm relies on the accuracy of the learners' approximation of the global solution and gradient.  In particular, in the commonly used strategy of single consensus averaging, every learner collects the information from immediate neighboring learners only. In a large network with a diverse set of learners, the knowledge of a small group of learners in a local neighborhood may not be sufficient to build an accurate approximation of the global loss, leading to decision update directions that deviate from the global minimizer. This issue can be resolved by applying consensus averaging for multiple times. It is important to have a proper number of consensus iterations to ensure that the errors due to the lack of access to global information are sufficiently small.  Yet, an increase in the consensus iterations leads to a larger communication overhead. Therefore, this requires a careful design to balance the need for information sharing and communication overhead. In DOMD-MADGC, we will show that a logarithmic number of consensus iterations in time suffices to achieve $O(1+C_T)$ dynamic regret.

The detailed procedure of DOMD-MADGC is as follows. At each round $t$, every learner node $i$ maintains a decision vector $x_{i,t}$, a local estimate of the global minimizer $y_{i,t}$, and a local approximate gradient of the global gradient $g_{i,t}$. At round $t$, every learner node $i$ updates $y_{i,t}$ using multiple consensus averaging iterations over the current local decisions. Let $W^{t}$ denote the weight matrix of the network used at round $t$. Let $z_{i,t}^{(0)} = x_{i,t}$ represent the initial message vector before beginning of consensus averaging. Let $z_{i,t}^{(K_t)}$ denote the vector output at learner node $i$ in round $t$, after $K_t$ consensus iterations. Then, at iteration $K_t$, after collecting messages from all its neighboring learners, learner $i$'s estimate of the global minimizer is
 \begin{align} \label{q_cons_y}
z_{i,t}^{(K_t)} = \sum_{j =1}^{n} (W^{t})_{ij} z_{j,t}^{(K_t-1)} = \sum_{j =1}^{n} ((W^{t})^{K_t})_{ij} z_{j,t}^{(0)},
\end{align}
 Every learner $i$ updates its approximation of the global minimizer after $K_t = \lceil \frac{-2 \log t}{ \log \sigma_2(W^t)}  \rceil$ consensus iterations, where $\sigma_2(W^t)$ is the second largest singular value of the communication network graph at round $t$, i.e., $W^t$.  Thus, we have $y_{i,t} = z_{i,t}^{(K_t)}$, which is represented by 
\begin{align} \label{q_cons_y_2}
    y_{i,t} = \sum_{j=1}^{n} ((W^{t})^{K_t})_{ij} x_{j,t}.
\end{align}

Then, the loss functions are revealed to the learners individually. Every learner $i$ computes the gradient of the local loss function $\nabla f_{t,i}(y_{t,i})$, which is evaluated at its estimate of the global minimizer $y_{i,t}$. Next, every learner $i$ updates its estimate of the global gradient using $K_t$ consensus averaging iterations over the local gradients, which leads to
\begin{align} \label{q_cons_g}
    g_{i,t} = \sum_{j=1}^{n} ((W^{t})^{K_t})_{ij} \nabla f_{j,t}(y_{j,t}).
\end{align}
The objective of this step is to enable the learners to access an estimate of the global gradient, which is collaboratively built based on the knowledge of the individual local learner as well as the neighboring learners.

After computing the local estimate vectors $y_{i,t}$ and $g_{i,t}$, learner node $i$ updates its local decision variable by
\begin{align}\label{q_update_main}
    x_{i, t+1} &= \mathbf{MD}_\eta(g_{i,t}, y_{i,t} ) \overset{\Delta}{=} \underset{{x \in \mathcal{X}}}{{\text{argmin}}}\Big\{\langle x, g_{i,t} \rangle + \frac{1}{\eta} D_r(x, y_{i,t}) \Big\},
\end{align}
where $\mathbf{MD}_\eta(g, y)$ represents the mirror descent update with a step size of $\eta$, and $D_{r}(x,y)$ is the Bregman divergence between $x$ and $y$ corresponding to the regularization function $r(x)$. Recall that $x_{i,t}$ and $y_{i,t}$ are respectively the local decision and the local estimate of the global minimizer at round~$t$. Thus, \eqref{q_update_main} suggests that every learner $i$ aims to stay close to the locally estimated global minimizer $y_{i,t}$ as measured by the Bregman divergence, while taking a step in a direction close to $g_{i,t}$ to reduce the local estimate of the global loss function at the current round. 

The above procedure of DOMD-MADGC is summarized in Algorithm~\ref{alg1}.

\begin{algorithm}[t] \caption{DOMD-MADGC} \label{alg1}
	\begin{algorithmic} 
		\STATE \textbf{Input:} Arbitrary initialization of $\{x_{i,1}\} \in \mathcal{X}$; step size $\eta$; time horizon $T$. 
		\STATE \textbf{Output:} Sequence of decisions $\{x_{i,t}, y_{i,t} : 1 \leq t \leq T\}$.
		\STATE 1: \textbf{for} $t = 1, 2, \ldots, T$ \textbf{do}\\
		\STATE 2: \quad set $K_t = \lceil \frac{-2 \log t}{ \log \sigma_2(W^t)}  \rceil$\\   
		\STATE 3: \quad $y_{i,t} = \sum_{j=1}^{n} ((W^{t})^{K_t})_{ij} x_{j,t}$\\   
		\STATE 4: \quad $g_{i,t} = \sum_{j=1}^{n} ((W^{t})^{K_t})_{ij} \nabla f_{j,t}(y_{j,t})$
		\STATE 5: \quad $x_{i, t+1} = \underset{{x \in \mathcal{X}}}{{\text{argmin}}}\Big\{\langle x, g_{i,t} \rangle + \frac{1}{\eta} D_r(x, y_{i,t}) \Big\}$
		\STATE 6: \textbf{end for}
	\end{algorithmic}
\end{algorithm}

\subsection{Improved Dynamic Regret}
\label{subsection_improved_regret}

In this section, we show how DOMD-MADGC improves the dynamic regret via multiple consensus averaging over both the local gradients and the local decisions. We first start by stating several standard assumptions commonly used in the literature of learning theory.

The following set of assumptions are commonly used in the literature after the group of studies began by \citep{hazan2007logarithmic, shalev2007logarithmic}, to provide stronger regret bounds by constraining the curvature of loss functions. For example, they are also used in \citep{chang2020unconstrained, hendrikx2020statistically, tarokh2020distributed}.

\textit{Assumption 1.} The loss functions $f_{i,t}(x)$ are $\lambda$-strongly convex and $\beta$-smooth. The loss functions have bounded gradients, i.e., $\|\nabla f_{i,t}(x)\| \leq G $. The regularization function $r(x)$ is $\mu$-strongly convex and $\mu'$-smooth.

The next set of assumptions pertain to the distributed nature of the learning network. Following previous studies on distributed optimization with time-varying network topology, we make several standard assumptions on topology graph $\mathcal{G}_t$. 

\textit{Assumption 2.} 
The weight matrix $W^t$ is doubly stochastic, i.e.,
\begin{align} \label{q_doubly_stochastic}
\sum_{j = 1}^{n} W^t_{ij} = \sum_{i = 1}^{n} W^t_{ij} = 1, \quad \forall t \geq 0.
\end{align}
Also, there exists a path from any learner $i$ to any learner $j$, i.e., the network is connected in each online round.

Assumption 2 ensures that the entries of the matrix  $(W^t)^K = \underbrace{ W^t \ldots W^t}_{K~\text{times}}$ are close to $1/n$, which is expressed by
\begin{align} \label{q_weight_matrix_1}
    \sum_{j = 1}^{n} \biggl| \bigg((W^t)^K\bigg)_{ij} - \frac{1}{n}  \biggl| \leq  \sqrt{n} \sigma_2^K(W^t),
\end{align}
where $\sigma_2(W^t)$ is the second largest singular value of the weight matrix $W^t$. This is a standard property of doubly stochastic matrices \citep{horn2012matrix}.

Now, we are ready to analyze the performance of DOMD-MADGC. We begin by bounding the distance between the exact average of the local decisions and the current global minimizer in the following lemma.

\begin{lemma} \label{lem_diff}
Under Assumptions 1 and 2, the sequence of decisions generated by DOMD-MADGC with the step size $\eta < \frac{\mu'}{\lambda}$ satisfy the following bound:
\begin{align} \label{lem_eq}
\| \bar{x}_{t+1} - x^*_t \| \leq \rho \| \bar{x}_{t} - x^*_t \| + \| \Delta_t \| + \| \delta_t \|, 
\end{align}
where $x^*_t$ is the minimizer of $f_t(x)$, $\bar{x}_t = \frac{1}{n} \sum_{i=1}^{n} x_{i, t}$ represents the exact average of local decisions, $\rho = \frac{\mu' - \eta \lambda}{\mu}$, and $\Delta_t$ and $\delta_t$ are given by
\begin{align} \label{q_network_errors}
& \Delta_t = \bar{x}_{t+1} - \mathbf{MD}_{\eta}(\bar{g}_t, \bar{x}_t), \\
& \delta_t = \mathbf{MD}_{\eta}(\bar{g}_t, \bar{x}_t) - \mathbf{MD}_{\eta}(\nabla f_t(\bar{x}_t), \bar{x}_t).
\end{align}
\end{lemma}

The following is a proof sketch of Lemma~\ref{lem_diff}. We first decompose the left-hand side of \eqref{lem_eq} into multiple terms, each representing a distinct source of error. In particular, the time-varying minimizers and limited local knowledge are two prominent source of errors. The error term due to the first source, called tracking error, measures the distance between the time-varying minimizer and the local decision of a learner, if it had access to exact average of the local decisions and global loss functions. Using the smooth duality property, from Lemma~2.19 in \citep{shalev2012online}, we bound the tracking error by $\| \bar{x}_{t} - x^*_t \|$. 
The other error terms are denoted by $\Delta_t$ and $\delta_t$. They are due to the second source of error,  which reflects the inadequacy of the knowledge of local learners. Thus, they are called network errors. Since the learners do not have access to the gradient of the global loss and the exact average of decisions $\bar{x}_t$, their approximation of these quantities will always have some network errors.   
The detailed  proof of Lemma~\ref{lem_diff} is given in App.~\ref{appx_lem_diff}. 

The result in Lemma~\ref{lem_diff} implies that after one round of DOMD-MADGC, the distance between the exact average of the local decisions $\bar{x}_{t+1}$ and the current global minimizer $x^*_t$ is less than a $\rho$-fraction of the distance between $\bar{x}_t$ and $x^*_t$, plus the network errors $\Delta_t$ and $\delta_t$. 

The following lemma provides an alternative form of update for mirror descent, which we require in our regret analysis.

\begin{lemma}\label{lemma_aux} Suppose $r(x)$ is strongly convex and $y$ satisfies $\nabla r(y) = \nabla r(u) - \eta \nabla f(l)$ for some convex function $f(x)$ and step size $\eta$. We have
\begin{align}
    \underset{x \in \mathcal{X}}{ \mathrm{argmin}}\bigg\{\langle \nabla f(l), x \rangle + \frac{1}{\eta} D_r(x, u) \bigg\} = \underset{x \in \mathcal{X}}{ \mathrm{argmin}}~D_r(x, y). \nonumber
\end{align}
\end{lemma}
Lemma~\ref{lemma_aux} is given as \textit{Proposition 17} in \citep{chiang12}.

We are now ready to upper bound the dynamic regret in Theorem~\ref{thm_reg}.

\begin{theorem} \label{thm_reg}
Under Assumptions 1 and 2, the dynamic regret of DOMD-MADGC, with $K_t = \lceil \frac{-2 \log t}{\log \sigma_2(W^t)} \rceil$ and a  fixed step size $ \frac{\mu' - \mu}{\lambda} < \eta < \frac{\mu'}{\lambda}$, satisfies 
\begin{align}
\mathrm{Reg}_T^\mathrm{d} & \leq    \frac{ G \| \bar{x}_1 - x^*_1  \| }{ 1 - \rho} +  \frac{G}{1-\rho} \sum_{t=2} ^{T}\| x^*_t - x^*_{t-1} \| 
 + \Big( \frac{GR \mu' + \eta G^2 + \eta \lambda G R }{\mu} + GR \Big) \frac{\sqrt{n} \pi^2}{6} \nonumber
\end{align}
where  $R = \max_{x \in \mathcal{X}} \|x\|$, and $C_T$ is the path length defined in \eqref{q_pl_def}.
\end{theorem}

The proof of Theorem~\ref{thm_reg} is given in App.~\ref{proof_thm}.

\textit{Remark 1.} Theorem~\ref{thm_reg} indicates that the dynamic regret of DOMD-MADGC is bounded by $O(1 + C_T)$, where $C_T$ is the path length, defined in \eqref{q_pl_def}.

\textit{Remark 2.} Recall that the dynamic regret of the distributed online mirror descent algorithm in \citep{shahrampour2017distributed} is upper bounded by $O(\sqrt{ T (1 + C_T)  })$. The $O(1+C_T)$ dynamic regret bound of DOMD-MADGC is smaller than $O(\sqrt{ T (1 + C_T)  })$ as long as $C_T$ is in the order of $o(T)$. Furthermore, in contrast to \citep{shahrampour2017distributed}, our regret bound can be achieved by setting a constant step size without any prior knowledge of the path length $C_T$ or the number of rounds $T$.  

\textit{Remark 3.} Recent studies on distributed online convex optimization have shown that the upper bound on the dynamic regret can be as tight as $O(1+ C_T + V_T)$ or $O(\log T (1 + C_T) )$ when the loss functions are strongly convex and smooth \citep{tarokh2020distributed, dixit2019distributed}. Theorem~~\ref{thm_reg} shows that by employing multiple consensus averaging iterations over both local decisions and local gradients, DOMD-MADGC can improve the dynamic regret bound to $O(1+ C_T)$. In particular, gradient exchange allows the learners to approximate and track the global gradient, which is boosted by multiple consensus iterations to allow them to collect the information of a larger number of nodes. These two steps together enable the learners to more precisely track the time-varying global minimizer, which lead to the improved regret bound of $O(1+ C_T)$.

\section{Discussion}

We now turn to investigate the effect of network topology and weight matrix on the communication load. Recall that $W^t$ represents the weight matrix of the network topology in online round $t$. The largest singular value $\sigma_1(W^t) = 1$ since $W^t$ is doubly stochastic. The number of consensus steps per round $K_t = \lceil \frac{-2 \log t}{ \log \sigma_2(W^t)}  \rceil \leq \lceil \frac{2 \log t}{ 1 - \sigma_2(W^t) } \rceil$, which follows from the fact that $\log \sigma^{-1}_2(W^t) \geq 1 - \sigma_2(W^t)$. Therefore, the parameter $K_t$ is controlled by the spectral gap of the weight matrix $W^t$, i.e., $\gamma(W^t) = 1 - \sigma_2(W^t)$.

Using the spectral gap, we can derive explicit bounds on the number of consensus steps for several interesting networks. 

\textit{Regular Grids:} The grid graph is formed by placing the nodes on a two dimensional plane. In $d$-connected regular grid, each node directly connects with its $d$ nearest neighbors in axis-aligned directions. The grid networks have often been used to model cluster computing and sensor networks.
For a $\sqrt{n}$-by-$\sqrt{n}$ $d$-regular grid network the second-largest singular value is bounded by $\sigma_2(W^t) = 1 - \Theta( \frac{k^2}{n} )$ \citep{duchi12dual, chung1997spectral}. In this case, the number of consensus steps per round is bounded by $K_t \in O( \frac{n \log t}{d^2} )$.

\textit{Random Geometric Graphs:} A random geometric graph can be constructed by placing the nodes uniformly on a two dimensional plane and connecting any pair of nodes whose distance is less than some pre-determined radius $r$. These graphs are commonly used to model the connectivity pattern of a set of wireless devices, such as in mobile communication networks. 
The properties of a random geometric graph with $r = (({\log}^{1+\epsilon} n) / n)^{1/2}$ and $\epsilon > 0$ are analyzed in \citep{von2014hitting}. In this case, the second-largest singular value of $W^t$ is bounded by $\sigma_2(W^t) = 1 - \Omega( \frac{{\log}^{1+\epsilon} n}{n})$, so $K_t$ is in the order of $O(\frac{n \log t}{\log n})$.

\textbf{Communication cost vs. regret bound:} The communication cost of DOMD-MADGC, which is defined as the overall number of consensus iterations, is $O(T \log T)$ after $T$ rounds. In comparison, prior works based on a single consensus iteration per online round has a communication cost of $O(T)$ after $T$ rounds. On the other hand, those prior works arrive at an upper bound of $O(\sqrt{T}(1+ C_T))$ on the dynamic regret when the prior knowledge of $C_T$ is not available, which has an additional $\sqrt{T}$ factor compared with our bound of $O(1+C_T)$. 
Thus, the extra communication cost of $O(\log T)$ can be seen as the price of removing the $O(\sqrt{T})$ factor in the dynamic regret bound. 

\section{Conclusion}

In this paper, we have presented a novel algorithm for distributed online optimization, with an aim to improve the dynamic regret. The proposed DOMD-MADGC runs multiple consensus averaging iterations over both the local decisions and the local gradients, which allows the distributed learners to accurately estimate the gradients of the time-varying global loss functions. Furthermore, our algorithm does not require any prior knowledge of the regularity measures, such as $C_T$, or the number of rounds $T$. Our theoretical analysis shows that the dynamic regret of DOMD-MADGC is bounded by $O(1+C_T)$, which is the best known result for distributed online optimization. 

\clearpage\newpage
\bibliographystyle{plainnat}
\bibliography{refs}

\begin{thebibliography}{41}
\providecommand{\natexlab}[1]{#1}
\providecommand{\url}[1]{\texttt{#1}}
\expandafter\ifx\csname urlstyle\endcsname\relax
  \providecommand{\doi}[1]{doi: #1}\else
  \providecommand{\doi}{doi: \begingroup \urlstyle{rm}\Url}\fi

\bibitem[Akbari et~al.(2015)Akbari, Gharesifard, and
  Linder]{akbari2015distributed}
Mohammad Akbari, Bahman Gharesifard, and Tam{\'a}s Linder.
\newblock Distributed online convex optimization on time-varying directed
  graphs.
\newblock \emph{IEEE Transactions on Control of Network Systems}, 4\penalty0
  (3):\penalty0 417--428, 2015.

\bibitem[Besbes et~al.(2015)Besbes, Gur, and Zeevi]{besbes2015non}
Omar Besbes, Yonatan Gur, and Assaf Zeevi.
\newblock Non-stationary stochastic optimization.
\newblock \emph{Operations Research}, 63\penalty0 (5):\penalty0 1227--1244,
  2015.

\bibitem[Boyd et~al.(2011)Boyd, Parikh, and Chu]{boyd2011distributed}
Stephen Boyd, Neal Parikh, and Eric Chu.
\newblock Distributed optimization and statistical learning via the alternating
  direction method of multipliers.
\newblock \emph{Foundations and Trends in Machine Learning}, 2011.

\bibitem[Campolongo and Orabona(2020)]{campolongo2020temporal}
Nicol\`{o} Campolongo and Francesco Orabona.
\newblock Temporal variability in implicit online learning.
\newblock In \emph{Advances in Neural Information Processing Systems}, 2020.

\bibitem[Campolongo and Orabona(2021)]{campolongo2021closer}
Nicol{\`o} Campolongo and Francesco Orabona.
\newblock A closer look at temporal variability in dynamic online learning.
\newblock \emph{arXiv preprint, arXiv:2102.07666}, 2021.

\bibitem[Cesa-Bianchi and Lugosi(2006)]{cesa2006prediction}
Nicolo Cesa-Bianchi and G{\'a}bor Lugosi.
\newblock \emph{Prediction, Learning, and Games}.
\newblock Cambridge university press, 2006.

\bibitem[Chang and Lin(2011)]{chang2011libsvm}
Chih-Chung Chang and Chih-Jen Lin.
\newblock {LIBSVM}: a library for support vector machines.
\newblock \emph{ACM Transactions on Intelligent Systems and Technology},
  2\penalty0 (3):\penalty0 1--27, 2011.

\bibitem[Chang and Shahrampour(2020)]{chang2020unconstrained}
Ting-Jui Chang and Shahin Shahrampour.
\newblock Unconstrained online optimization: Dynamic regret analysis of
  strongly convex and smooth problems.
\newblock \emph{arXiv preprint, arXiv:2006.03912}, 2020.

\bibitem[Chiang et~al.(2012)Chiang, Yang, Lee, Mahdavi, Lu, Jin, and
  Zhu]{chiang12}
Chao-Kai Chiang, Tianbao Yang, Chia-Jung Lee, Mehrdad Mahdavi, Chi-Jen Lu, Rong
  Jin, and Shenghuo Zhu.
\newblock Online optimization with gradual variations.
\newblock In \emph{Proceedings of Annual Conference on Learning Theory}, 2012.

\bibitem[Chung and Graham(1997)]{chung1997spectral}
Fan~RK Chung and Fan~Chung Graham.
\newblock \emph{Spectral graph theory}.
\newblock Number~92. American Mathematical Soc., 1997.

\bibitem[Dixit et~al.(2019)Dixit, Bedi, Rajawat, and
  Koppel]{dixit2019distributed}
Rishabh Dixit, Amrit~Singh Bedi, Ketan Rajawat, and Alec Koppel.
\newblock Distributed online learning over time-varying graphs via proximal
  gradient descent.
\newblock In \emph{Proceedings of IEEE International Conference on Decision and
  Control}, 2019.

\bibitem[Duchi et~al.(2010)Duchi, Shalev-Shwartz, Singer, and
  Tewari]{duchi2010composite}
John~C Duchi, Shai Shalev-Shwartz, Yoram Singer, and Ambuj Tewari.
\newblock Composite objective mirror descent.
\newblock In \emph{Proceedings of the Conference on Learning Theory}, 2010.

\bibitem[Duchi et~al.(2012)Duchi, Agarwal, and Wainwright]{duchi12dual}
John~C. Duchi, Alekh Agarwal, and Martin~J. Wainwright.
\newblock Dual averaging for distributed optimization: Convergence analysis and
  network scaling.
\newblock \emph{IEEE Transactions on Automatic Control}, 57\penalty0 (3), 2012.

\bibitem[Hall and Willett(2013)]{hall2013dynamical}
Eric Hall and Rebecca Willett.
\newblock Dynamical models and tracking regret in online convex programming.
\newblock In \emph{Proceedings of International Conference on Machine
  Learning}, 2013.

\bibitem[Hall and Willett(2015)]{hall2015}
Eric~C Hall and Rebecca~M Willett.
\newblock Online convex optimization in dynamic environments.
\newblock \emph{IEEE Journal of Selected Topics in Signal Processing},
  9\penalty0 (4):\penalty0 647--662, 2015.

\bibitem[Hazan(2019)]{hazan2019introduction}
Elad Hazan.
\newblock Introduction to online convex optimization.
\newblock \emph{arXiv preprint, arXiv:1909.05207}, 2019.

\bibitem[Hazan et~al.(2006)Hazan, Kalai, Kale, and
  Agarwal]{hazan2007logarithmic}
Elad Hazan, Adam~Tauman Kalai, Satyen Kale, and Amit Agarwal.
\newblock Logarithmic regret algorithms for online convex optimization.
\newblock In \emph{Proceedings of the Conference on Learning Theory}, 2006.

\bibitem[Hendrikx et~al.(2020)Hendrikx, Xiao, Bubeck, Bach, and
  Massoulie]{hendrikx2020statistically}
Hadrien Hendrikx, Lin Xiao, Sebastien Bubeck, Francis Bach, and Laurent
  Massoulie.
\newblock Statistically preconditioned accelerated gradient method for
  distributed optimization.
\newblock 2020.

\bibitem[Horn and Johnson(2012)]{horn2012matrix}
Roger~A Horn and Charles~R Johnson.
\newblock \emph{Matrix analysis}.
\newblock Cambridge university press, 2012.

\bibitem[Hosseini et~al.(2013)Hosseini, Chapman, and
  Mesbahi]{hosseini2013online}
Saghar Hosseini, Airlie Chapman, and Mehran Mesbahi.
\newblock Online distributed optimization via dual averaging.
\newblock In \emph{Proceedings of IEEE International Conference on Decision and
  Control}, 2013.

\bibitem[Jadbabaie et~al.(2015)Jadbabaie, Rakhlin, Shahrampour, and
  Sridharan]{jadbabaie2015online}
Ali Jadbabaie, Alexander Rakhlin, Shahin Shahrampour, and Karthik Sridharan.
\newblock {Online optimization:Competing with dynamic comparators}.
\newblock In \emph{Proceedings of the International Conference on Artificial
  Intelligence and Statistics}, 2015.

\bibitem[Krizhevsky(2009)]{krizhevsky2009}
Alex Krizhevsky.
\newblock Learning multiple layers of features from tiny images.
\newblock 2009.
\newblock URL
  \url{http://www.cs.toronto.edu/~kriz/learning-features-2009-TR.pdf}.

\bibitem[Lee et~al.(2016)Lee, Nedi{\'c}, and Raginsky]{lee2016coordinate}
Soomin Lee, Angelia Nedi{\'c}, and Maxim Raginsky.
\newblock Coordinate dual averaging for decentralized online optimization with
  nonseparable global objectives.
\newblock \emph{IEEE Transactions on Control of Network Systems}, 5\penalty0
  (1):\penalty0 34--44, 2016.

\bibitem[Li et~al.(2020)Li, Yi, and Xie]{li2020distributed}
Xiuxian Li, Xinlei Yi, and Lihua Xie.
\newblock Distributed online convex optimization with an aggregative variable.
\newblock \emph{arXiv preprint, arXiv:2007.06844}, 2020.

\bibitem[Lu et~al.(2019)Lu, Jing, and Wang]{lu2019online}
Kaihong Lu, Gangshan Jing, and Long Wang.
\newblock Online distributed optimization with strongly pseudoconvex-sum cost
  functions.
\newblock \emph{IEEE Transactions on Automatic Control}, 65\penalty0
  (1):\penalty0 426--433, 2019.

\bibitem[Mateos-N{\'u}nez and Cort{\'e}s(2014)]{mateos2014distributed}
David Mateos-N{\'u}nez and Jorge Cort{\'e}s.
\newblock Distributed online convex optimization over jointly connected
  digraphs.
\newblock \emph{IEEE Transactions on Network Science and Engineering},
  1\penalty0 (1):\penalty0 23--37, 2014.

\bibitem[Mokhtari et~al.(2016)Mokhtari, Shahrampour, Jadbabaie, and
  Ribeiro]{mokhtari2016online}
Aryan Mokhtari, Shahin Shahrampour, Ali Jadbabaie, and Alejandro Ribeiro.
\newblock Online optimization in dynamic environments: Improved regret rates
  for strongly convex problems.
\newblock In \emph{Proceedings of IEEE International Conference on Decision and
  Control}, 2016.

\bibitem[Orabona(2019)]{orabona2019modern}
Francesco Orabona.
\newblock A modern introduction to online learning.
\newblock \emph{arXiv preprint, arXiv:1912.13213}, 2019.

\bibitem[Rabbat(2015)]{rabbat15}
Michael Rabbat.
\newblock Multi-agent mirror descent for decentralized stochastic optimization.
\newblock In \emph{IEEE International Workshop on Computational Advances in
  Multi-Sensor Adaptive Processing}, 2015.

\bibitem[Shahrampour and Jadbabaie(2018)]{shahrampour2017distributed}
Shahin Shahrampour and Ali Jadbabaie.
\newblock Distributed online optimization in dynamic environments using mirror
  descent.
\newblock \emph{IEEE Transactions on Automatic Control}, 63\penalty0
  (3):\penalty0 714--725, 2018.

\bibitem[Shalev-Shwartz(2012)]{shalev2012online}
Shai Shalev-Shwartz.
\newblock Online learning and online convex optimization.
\newblock \emph{Foundations and Trends in Machine Learning}, 4\penalty0
  (2):\penalty0 107--194, 2012.

\bibitem[Shalev-Shwartz and Singer(2007)]{shalev2007logarithmic}
Shai Shalev-Shwartz and Yoram Singer.
\newblock \emph{Logarithmic Regret Algorithms for Strongly Convex Repeated
  Games}.
\newblock The Hebrew University, 2007.

\bibitem[Tsianos and Rabbat(2016)]{tsianos2016efficient}
Konstantinos~I Tsianos and Michael~G Rabbat.
\newblock Efficient distributed online prediction and stochastic optimization
  with approximate distributed averaging.
\newblock \emph{IEEE Transactions on Signal and Information Processing over
  Networks}, 2\penalty0 (4):\penalty0 489--506, 2016.

\bibitem[Von~Luxburg et~al.(2014)Von~Luxburg, Radl, and Hein]{von2014hitting}
Ulrike Von~Luxburg, Agnes Radl, and Matthias Hein.
\newblock Hitting and commute times in large random neighborhood graphs.
\newblock \emph{The Journal of Machine Learning Research}, 15\penalty0
  (1):\penalty0 1751--1798, 2014.

\bibitem[Xi et~al.(2014)Xi, Wu, and Khan]{xi5526distributed}
C~Xi, Q~Wu, and UA~Khan.
\newblock Distributed mirror descent over directed graphs.
\newblock \emph{arXiv preprint, arXiv:1412.5526}, 2014.

\bibitem[Xiao(2010)]{xiao2010dual}
Lin Xiao.
\newblock Dual averaging methods for regularized stochastic learning and online
  optimization.
\newblock \emph{Journal of Machine Learning Research}, 11\penalty0
  (88):\penalty0 2543--2596, 2010.

\bibitem[Yang et~al.(2016)Yang, Zhang, Jin, and Yi]{yangb16}
Tianbao Yang, Lijun Zhang, Rong Jin, and Jinfeng Yi.
\newblock Tracking slowly moving clairvoyant: Optimal dynamic regret of online
  learning with true and noisy gradient.
\newblock In \emph{Proceedings of International Conference on Machine
  Learning}, 2016.

\bibitem[Yuan et~al.(2020)Yuan, Hong, Ho, and Xu]{yuan2020distributed}
Deming Yuan, Yiguang Hong, Daniel~WC Ho, and Shengyuan Xu.
\newblock Distributed mirror descent for online composite optimization.
\newblock \emph{IEEE Transactions on Automatic Control}, 2020.

\bibitem[Zhang et~al.(2017)Zhang, Yang, Yi, Jin, and Zhou]{zhang2017improved}
Lijun Zhang, Tianbao Yang, Jinfeng Yi, Rong Jin, and Zhi-Hua Zhou.
\newblock Improved dynamic regret for non-degenerate functions.
\newblock In \emph{Proceedings of the International Conference on Advances in
  Neural Information Processing Systems}, 2017.

\bibitem[Zhang et~al.(2020)Zhang, Ravier, Tarokh, and
  Zavlanos]{tarokh2020distributed}
Yan Zhang, Robert~J. Ravier, Vahid Tarokh, and Michael~M. Zavlanos.
\newblock Distributed online convex optimization with improved dynamic regret.
\newblock \emph{arXiv preprint, arXiv:1911.05127}, 2020.

\bibitem[Zinkevich(2003)]{zinkevich2003online}
Martin Zinkevich.
\newblock Online convex programming and generalized infinitesimal gradient
  ascent.
\newblock In \emph{Proceedings of the International Conference on Machine
  Learning}, 2003.

\end{thebibliography}

\clearpage\newpage
\appendix

\section{Proof of Lemma~\ref{lem_diff}}\label{appx_lem_diff}
\allowdisplaybreaks
We begin by considering the difference between $\bar{x}_{t+1} = \frac{1}{n} \sum_{i=1}^{n} x_{i, t+1}$ and $x^*_t$.
Using the update in \eqref{q_update_main}, i.e., $x_{i, t+1} = \mathbf{MD}_\eta(g_{i,t}, y_{i,t} )$, and by adding and subtracting several terms we obtain
\begin{align}
    \| \bar{x}_{t+1} - x^*_t  \| & = \Big\| \frac{1}{n} \sum_{i=1}^{n} \mathbf{MD}_{\eta}(g_{i,t}, y_{i,t}) - \mathbf{MD}_{\eta}( \nabla f_t(x^*_t), x^*_t)   \Big\| \nonumber\\
    & \leq \Big\| \frac{1}{n} \sum_{i=1}^{n} \mathbf{MD}_{\eta}(g_{i,t}, y_{i,t}) - \mathbf{MD}_{\eta}( \bar{g}_t, \bar{x}_t) + \mathbf{MD}_{\eta}( \bar{g}_t, \bar{x}_t) -\mathbf{MD}_{\eta}( \nabla f_t(x^*_t), x^*_t)   \Big\| \nonumber\\ 
    & \leq \Big \|\mathbf{MD}_{\eta}( \bar{g}_t, \bar{x}_t) - \mathbf{MD}_{\eta}( \nabla f_t(\bar{x}_t), \bar{x}_t) + \mathbf{MD}_{\eta}( \nabla f_t(\bar{x}_t), \bar{x}_t) - \mathbf{MD}_{\eta}( \nabla f_t(x^*_t), x^*_t)   \Big\| \nonumber\\
    & \quad\quad + \Big\| \frac{1}{n} \sum_{i=1}^{n} \mathbf{MD}_{\eta}(g_{i,t}, y_{i,t}) - \mathbf{MD}_{\eta}( \bar{g}_t, \bar{x}_t) \Big\|  \nonumber \\
    & \leq    \Big\| \mathbf{MD}_{\eta}( \nabla f_t(\bar{x}_t), \bar{x}_t) - \mathbf{MD}_{\eta}( \nabla f_t(x^*_t), x^*_t)   \Big\| + \Big\|\mathbf{MD}_{\eta}( \bar{g}_t, \bar{x}_t) - \mathbf{MD}_{\eta}( \nabla f_t(\bar{x}_t), \bar{x}_t) \Big\| \nonumber \\
    & \quad\quad + \Big\| \frac{1}{n} \sum_{i=1}^{n} \mathbf{MD}_{\eta}(g_{i,t}, y_{i,t}) - \mathbf{MD}_{\eta}( \bar{g}_t, \bar{x}_t) \Big\|, \label{aaa}
\end{align}
where the first line follows from the fact that $x^*_t = \mathbf{MD}_{\eta}( \nabla f_t(x^*_t), x^*_t)$. We bound each of the terms on the right hand-side of \eqref{aaa} separately. The first term of the above inequality can be expanded as
\begin{align}
      \Big\| \mathbf{MD}_{\eta}( \nabla f_t(\bar{x}_t), \bar{x}_t) & - \mathbf{MD}_{\eta}( \nabla f_t(x^*_t), x^*_t)  \Big\|^2 \nonumber\\
       & \leq 
      \| \nabla r^*( \nabla r (\bar{x}_t) - \eta \nabla f_t(\bar{x}_t) ) - \nabla r^*( \nabla r (x^*_t) - \eta \nabla f_t(x^*_t) )  \|^2 \nonumber\\
      & \leq \frac{1}{\mu^2} \| \nabla r (\bar{x}_t) - \eta \nabla f_t(\bar{x}_t) - (\nabla r (x^*_t) - \eta \nabla f_t(x^*_t)) \|^2,
 \label{q_lem_2}
\end{align}
where the first line is obtained using the alternate form of mirror descent update stated in Lemma 2, and the fact that the inverse of $\nabla r(x)$ is $\nabla r^*(x)$ when $r(x)$ is strongly convex. Here $r^*(x)$ denotes the conjugate function \citep{shalev2012online}. We have also used the fact that the distance of the projection of two points into a convex set is smaller than the distance between these unprojected points. Note that since $r(x)$ is $\mu$-strongly convex, its conjugate $r^*(x)$ is $1/\mu$-smooth \citep{shalev2012online}.

We now proceed to bound the right hand-side of \eqref{q_lem_2}. The smoothness of $r(x)$ implies
\begin{align}
    r(y) \leq r(x) + \langle \nabla r(x), y - x \rangle + \frac{\mu'}{2} \|y-x\|^2, \forall x, y \in \mathcal{X}.  \label{q_lem_3}
\end{align}
In addition, $f_t(x)$ is $\lambda$-strongly convex, i.e., 
\begin{align}
    f_t(y) \geq f_t(x) + \langle \nabla f_t(x), y - x  \rangle + \frac{\lambda}{2} \| y-x \|^2, \forall x, y \in \mathcal{X}.  \label{q_lem_4}
\end{align}
We multiply \eqref{q_lem_4} by $-\eta$ and add it to \eqref{q_lem_3} to obtain
\begin{align}
    r(y) - \eta f_t(y) \leq r(x) - \eta f_t(x) + \langle \nabla r(x) - \eta \nabla f_t(x), y - x \rangle + \frac{\mu' - \eta \lambda}{2} \|y-x\|^2, \forall x, y \in \mathcal{X},  \label{q_lem_5}
\end{align}
which shows that the function $r(x) - \eta f_t(x)$ is $(\mu' - \eta\lambda)$-smooth when $\eta < \mu'/\lambda$. Therefore, we have
\begin{align}
   \| \nabla r (\bar{x}_t) - \eta \nabla f_t(\bar{x}_t) - (\nabla r (x^*_t) - \eta \nabla f_t(x^*_t)) \|^2 \leq (\mu' - \eta \lambda)^2 \| \bar{x}_t - x^*_t\|^2. \label{q_lem_6}
\end{align}
By combining the above inequality and \eqref{q_lem_2}, we have
\begin{align}
    \Big\| \mathbf{MD}_{\eta}( \nabla f_t(\bar{x}_t), \bar{x}_t) & - \mathbf{MD}_{\eta}( \nabla f_t(x^*_t), x^*_t)  \Big\|^2 \leq \rho^2 \| \bar{x}_t - x^*_t\|^2,  \label{q_lem_7}
\end{align}
where $\rho = \frac{\mu' - \eta \lambda}{\mu}$.
Substituting \eqref{q_lem_7} into \eqref{aaa}, and using the definitions of $ \Delta_t$ and $ \delta_t $ in \eqref{q_network_errors}, we get
\begin{align}
\| \bar{x}_{t+1} - x^*_t \| \leq \rho \| \bar{x}_{t} - x^*_t \| + \| \Delta_t \| + \| \delta_t \|.  
\end{align}


\section{Proof of Theorem~\ref{thm_reg}}\label{proof_thm}

We provide a proof sketch for Theorem~\ref{thm_reg}. We begin the proof by first bounding the network error term $\Delta_t$. Using the definition of $\Delta_t$ and simple norm properties, we have
\begin{align}
 & \| \Delta_t \| = \Big\| \frac{1}{n} \sum_{i=1}^{n} \mathbf{MD}_{\eta}(g_{i,t}, y_{i,t}) - \mathbf{MD}_{\eta}( \bar{g}_t, \bar{x}_t) \Big\| \nonumber\\
& \leq \frac{1}{n} \sum_{i=1}^{n} \Big\| \mathbf{MD}_{\eta}(g_{i,t}, y_{i,t}) - \mathbf{MD}_{\eta}( \bar{g}_t, \bar{x}_t) \Big\| \nonumber\\
& \leq \frac{1}{n} \sum_{i=1}^{n} \| \nabla r^*( \nabla r( y_{i,t} ) - \eta g_{i,t} ) - \nabla r^*( \nabla r( \bar{x}_t ) - \eta \bar{g}_t ) \|, \label{q_ax2_1} \raisetag{1.3cm}
\end{align}
where the last line follows the result of Lemma~\ref{lemma_aux}, and the fact that the inverse of $\nabla r(x)$ is $\nabla r^*(x)$ when $r(x)$ is strongly convex \citep{shalev2012online}. We note that $r^*(x)$ denotes the conjugate of $r(x)$, which is defined by $r^*(x) = \max_u \{ \langle u, x \rangle - r(u)\}$. We have also used the fact that the distance between the projection of two points into a convex set is less than that of unprojected points. In addition, since $r(x)$ is $\mu$-strongly convex, its conjugate $r^*(x)$ is $1/\mu$-smooth \citep{shalev2012online}. Thus, we have 
\begin{align}
     \| \nabla r^*( \nabla r( & y_{i,t} )  - \eta g_{i,t} ) - \nabla r^*( \nabla r( \bar{x}_t ) -  \eta \bar{g}_t ) \| \nonumber\\
    & \leq  \frac{1}{\mu} \| \nabla r( y_{i,t} ) - \eta g_{i,t}  -  \nabla r( \bar{x}_t ) + \eta \bar{g}_t  \| \nonumber\\
    & \leq  \frac{1}{\mu} \Big[ \| \nabla r( y_{i,t} ) -  \nabla r( \bar{x}_t ) \| + \eta \| g_{i,t} - \bar{g}_t \| \Big] \nonumber\\
    & \leq  \frac{\mu'}{\mu} \|  y_{i,t} - \bar{x}_t \| +   \frac{\eta}{\mu} \| g_{i,t} - \bar{g}_t \|, \label{q_ax2_1_2}
\end{align}
where the last line is obtained due to the smoothness of $r(x)$. 
The above inequality implies that the norm of the network error term $\Delta_t$, and ultimately the dynamic regret, are bounded by $\|  y_{i,t} - \bar{x}_t \|$, which indicates the need for consensus over local decisions, and by $\| g_{i,t} - \bar{g}_t \|$, which indicates the necessity of consensus over local gradients. The resultant error bound based on these terms shows the importance of having accurate estimation of the exact average of decisions and gradients, which motivates running a sufficient number of consensus iterations on both local decisions and local gradients to reduce the above terms. Thus, we combine \eqref{q_ax2_1} and \eqref{q_ax2_1_2}, and noting the number of consensus iterations $K_t = \lceil \frac{-2 \log t}{\log \sigma_2(W^t)} \rceil$,  to obtain 
\begin{align}    
  \| \Delta_t \|  & \leq \frac{\mu'}{n \mu} \sum_{i=1}^{n} \Big\| \sum_{j=1}^{n} \bigg( \bigg((W^t)^{K_t}\bigg)_{ij} - \frac{1}{n} \bigg) x_{j,t} \Big\| \nonumber\\
   &~~+ \frac{\eta}{n \mu } \sum_{i=1}^{n} \Big\| \sum_{j=1}^{n} \bigg( \bigg((W^t)^{K_t}\bigg)_{ij} - \frac{1}{n} \bigg) \nabla f_{j,t}(y_{j,t}) \Big\| \nonumber \\
   & \leq \frac{\mu'}{n \mu} \sum_{i=1}^{n}\sum_{j=1}^{n} \biggl| \bigg((W^t)^{K_t}\bigg)_{ij} - \frac{1}{n} \biggl| \Big\|x_{j,t} \Big\| \nonumber \\
   &~~+ \frac{\eta}{n \mu} \sum_{i=1}^{n}\sum_{j=1}^{n} \biggl| \bigg((W^t)^{K_t}\bigg)_{ij} - \frac{1}{n} \biggl|  \Big\| \nabla f_{j,t}(y_{j,t}) \Big\| \nonumber \\
   & \leq \frac{ \mu'}{ n\mu} \sum_{i=1}^{n} \sqrt{n} \sigma_2^{K_t}(W^t) R +  \frac{\eta}{n \mu} \sum_{i=1}^n \sqrt{n} \sigma_2^{K_t}(W^t) G \nonumber\\
   & \leq \frac{ \mu'}{\mu} \frac{ R \sqrt{n}}{t^2} +  \frac{\eta}{\mu} \frac{ G \sqrt{n}}{t^2}, \label{q_ax2_1_3}
\end{align}
where the third inequality follows from \eqref{q_weight_matrix_1}. The above inequality implies that the network error $\Delta_t$ is bounded by $O(\frac{1}{t^2})$. This is an important result, since it implies that $\Delta_t$ is summable over time.

In the next step, we bound the other error term $\delta_t$. Similar to the case with $\Delta_t$, we use an alternative form of mirror descent update, and the properties of the conjugate function $r^*(x)$. We show that the norm of  $\delta_t$ is bounded by $\|  y_{i,t} - \bar{x}_t \|$. Our analysis further shows that a logarithmically increasing number of consensus averaging iterations $K_t = \lceil \frac{-2 \log t}{\log \sigma_2(W^t)} \rceil$ is sufficient to bound $\delta_t$ also by $O(1/t^2)$. 
Next, using the result of Lemma~\ref{lem_diff}, and the upper bounds derived on $\Delta_t$ and $\delta_t$, we bound the cumulative distance between the exact average of decisions and the time-varying minimizer $\sum_{t=1}^{T} \| \bar{x}_{t} - x^*_t \|$. We choose a step size  $ \frac{\mu' - \mu}{\lambda} < \eta < \frac{\mu'}{\lambda}$ to ensure that $\rho < 1$ in order to bound $\sum_{t=1}^{T} \| \bar{x}_{t} - x^*_t \|$ by $O(1+C_T)$.

We note that $\Delta_t$ is bounded by $O(1/t^2)$, shown in \eqref{q_ax2_1_3}. Now we proceed to bound the other term on network error $\delta_t$ as follows:
\begin{align} 
\| \delta_t \| & = \| \mathbf{MD}_{\eta}( \bar{g}_t, \bar{x}_t) - \mathbf{MD}_{\eta}( \nabla f_t(\bar{x}_t), \bar{x}_t) \| \nonumber\\
    & \leq \| \nabla r^*( \nabla r(\bar{x}_t) - \eta \bar{g}_t ) - \nabla r^*( \nabla r(\bar{x}_t) - \eta \nabla f_t(\bar{x}_t) ) \| \nonumber\\
    & \leq \frac{1}{\mu} \|  \nabla r(\bar{x}_t) - \eta \bar{g}_t - \nabla r(\bar{x}_t) + \eta \nabla f_t(\bar{x}_t)  \| \nonumber \\
    & \leq \frac{\eta}{\mu}  \|  \bar{g}_t - \nabla f_t(\bar{x}_t)  \| \nonumber\\
    & \leq \frac{\eta}{\mu}  \Big\|  \frac{1}{n} \sum_{i=1}^{n} g_{i,t} - \frac{1}{n} \sum_{j=1}^{n} \nabla f_{j,t}(\bar{x}_t) \Big\|\label{q_ax2_2}
\end{align}
where the first inequality follows from the alternative form of mirror descent update in Lemma~\ref{lemma_aux}. We also note that since $r(x)$ is $\mu$-strongly convex, which implies its conjugate $r^*(x)$ is $1/\mu$-smooth. The last line shows that the network error $\delta_t$ relates to the distance between the average of the gradients of local learners and the gradient of the global loss function evaluated at $\bar{x}_t$. Therefore, we have
\begin{align}
 \| \delta_t \| & \leq \frac{\eta}{\mu} \Big\| \frac{1}{n} \sum_{i=1}^{n} \sum_{j=1}^{n} \bigg( (W^t)^{K_t} \bigg)_{ij}  \nabla f_{j,t}(y_{j,t})   -  \frac{1}{n} \sum_{j=1}^{n} \nabla f_{j,t}(\bar{x}_t) \Big\| \nonumber\\
    & \leq \frac{\eta}{\mu} \Big\| \frac{1}{n} \sum_{j=1}^{n} \nabla f_{j,t}(y_{j,t}) - \frac{1}{n} \sum_{j=1}^{n} \nabla f_{j,t}(\bar{x}_t) \Big\| \nonumber\\
    & \leq \frac{\eta}{\mu } \frac{1}{n} \sum_{j=1}^{n} \Big\| \nabla f_{j,t}(y_{j,t}) - \nabla f_{j,t}(\bar{x}_t) \Big\| \nonumber \\
    & \leq \frac{\eta}{\mu n}  \sum_{j=1}^{n} \lambda \| y_{j,t} - \bar{x}_t \|, \label{q_ax2_3}
\end{align}
where the second line follows from the fact that the matrix $(W^t)^{K_t} = \underbrace{ W^t \ldots W^t}_{K_t~\text{times}}$ is doubly stochastic. We have also used the strong convexity of local loss functions to obtain the right-hand side of \eqref{q_ax2_3}. The above inequality illustrates that the error $\delta_t$ is directly bounded by the distance between the local inexact average $y_{j,t}$ and the exact average of decisions $\bar{x}_t$. Thus, with $K_t$ consensus iterations, we have
\begin{align}
  \| \delta_t \|  & \leq \frac{\eta}{\mu n}  \sum_{j=1}^{n} \lambda \sum_{i=1}^{n} \biggl| \bigg((W^t)^{K_t}\bigg)_{ij} - \frac{1}{n} \biggl| \Big\|  x_{i,t} \Big\| \nonumber\\
    & \leq \frac{\eta}{\mu } \sqrt{n}  \lambda   \sigma_2^{K_t}(W^t) R \nonumber\\
    & \leq \frac{\eta \lambda  \sqrt{n} R }{\mu t^2}, \label{q_ax2_4}
\end{align}
where the last line is obtained using $K_t = \lceil \frac{-2 \log t}{\log \sigma_2(W^t)} \rceil$.
From \eqref{q_ax2_1_3} and \eqref{q_ax2_4}, we observe that both error terms $\Delta_t$ and $\delta_t$ are bounded by $O(1/t^2)$.

Next, we upper bound the cumulative distance between $\bar{x}_t$ and $x^*_t$. We first expand the summation and add and subtract several terms to obtain
\begin{align} \label{q_ax2_3_2}
\sum_{t=1}^{T} \| \bar{x}_{t} - x^*_t \| & = \| \bar{x}_1 - x^*_1  \| + \sum_{t=2}^{T} \| \bar{x}_t - x^*_{t-1} + x^*_{t-1}-x^*_t \| \nonumber\\
    & \leq \| \bar{x}_1 - x^*_1  \| + \sum_{t=2}^{T} \| \bar{x}_t - x^*_{t-1} \| + \sum_{t=2} ^{T}\| x^*_t - x^*_{t-1} \| \nonumber \\
    & \leq \| \bar{x}_1 - x^*_1  \| + \sum_{t=2}^{T} \Big[ \rho \| \bar{x}_{t-1} - x^*_{t-1}  \| + \| \Delta_t \| + \| \delta_t \| \Big] + \sum_{t=2} ^{T}\| x^*_t - x^*_{t-1} \| \nonumber\\
    & \leq \| \bar{x}_1 - x^*_1  \| + \rho \sum_{t=1}^{T} \| \bar{x}_{t} - x^*_{t}  \| +  \sum_{t=1}^{T} \Big( \frac{R \mu' + \eta G + \eta \lambda R }{\mu}  \Big) \frac{\sqrt{n}}{t^2} 
    + \sum_{t=2} ^{T}\| x^*_t - x^*_{t-1} \|,
\end{align}
where we have applied Lemma~\ref{lem_diff} in the third line. We have further used the network error bounds in \eqref{q_ax2_1_3} and \eqref{q_ax2_4} to obtain the right hand-side of \eqref{q_ax2_3_2}. The step size is chosen $\frac{\mu' - \mu}{\lambda} < \eta < \frac{\mu'}{\lambda} $ such that $\rho < 1$. Thus, we have
\begin{align}\label{q_ax2_4_2}
\sum_{t=1}^{T} \| \bar{x}_{t} - x^*_t \| \leq \frac{ \| \bar{x}_1 - x^*_1  \| }{ 1 - \rho} + \Big( \frac{R \mu' + \eta G + \eta \lambda R }{\mu}  \Big) \frac{ \sqrt{n} \pi^2}{6 (1 - \rho)} + \frac{1}{1-\rho} \sum_{t=2} ^{T}\| x^*_t - x^*_{t-1} \|.
\end{align}

Now we are ready to establish an upper bound on the dynamic regret. Since each learner maintains two as follows:
\begin{align}
\frac{1}{n} \sum_{i=1}^{n} \sum_{t=1}^T f_{t}(y_{i,t}) - \sum_{t=1}^T f_t(x^*_t) & = \frac{1}{n^2} \sum_{i=1}^{n}\sum_{j=1}^{n} \sum_{t=1}^T f_{j,t}(y_{i,t}) - \frac{1}{n} \sum_{t=1}^T\sum_{j=1}^n f_{j,t}(x^*_t)  \nonumber\\
& \leq \frac{1}{n^2} \sum_{i=1}^{n}\sum_{j=1}^{n} \sum_{t=1}^T f_{j,t}(y_{i,t}) - f_{j,t}(x^*_t) \nonumber \\
& \leq \frac{1}{n^2} \sum_{i=1}^{n}\sum_{j=1}^{n} \sum_{t=1}^T G \| y_{i,t} - x^*_t \| \nonumber \\
& \leq \frac{1}{n} \sum_{i=1}^{n} \sum_{t=1}^T G \| y_{i,t} - \bar{x}_t \| +  \sum_{t=1}^T G \| \bar{x}_t - x^*_t \|,  \label{q_ax2_5}
\end{align}
where we have used the Lipschitz continuity of local loss functions in the third line. It can be observed that the dynamic regret is a sum of two components. The first component reflects how close the local inexact average estimate is to the exact average of decisions. The second term measures the difference between the average of local decisions and the time-varying global minimizer, which is bounded in \eqref{q_ax2_4_2}. Therefore, we now only need to bound the first term. Using \eqref{q_weight_matrix_1}, we have  
\begin{align}
\frac{1}{n} \sum_{i=1}^{n} \sum_{t=1}^T G \| y_{i,t} - \bar{x}_t \| & \leq \frac{1}{n} \sum_{i=1}^{n} \sum_{t=1}^T G \sum_{j=1}^n \biggl| \bigg((W^t)^{K_t}\bigg)_{ij} - \frac{1}{n} \biggl| \Big\|x_{j,t} \Big\| \nonumber\\
& \leq  \sum_{t=1}^T G R \sqrt{n}  \sigma_2(W^t)^{K_t}  \nonumber\\
& \leq  G R \sqrt{n}  \sum_{t=1}^T \frac{1}{t^2}  \nonumber\\
& \leq  G R \sqrt{n}  \frac{\pi^2}{6}. \label{q_ax2_6}
\end{align}
Substituting \eqref{q_ax2_4_2} and \eqref{q_ax2_6} into \eqref{q_ax2_5} leads to the dynamic regret bound of
\begin{align}
\mathrm{Reg}_T^\mathrm{d}  \leq & G R \sqrt{n}  \frac{\pi^2}{6} +  G\frac{ \| \bar{x}_1 - x^*_1  \| }{ 1 - \rho}  \nonumber\\
&+  \Big( \frac{GR \mu' + \eta G^2 + \eta \lambda G R }{\mu}  \Big) \frac{ \sqrt{n} \pi^2}{6 (1 - \rho)} + \frac{G}{1-\rho} \sum_{t=2} ^{T}\| x^*_t - x^*_{t-1} \|.
\end{align}

\section{Experiments}
\label{sec_experiment}

We illustrate the performance of DOMD-MADGC with a series of numerical experiments on real datasets. 
In particular, we consider the problem of multi-class classification in the distributed online learning setting over  notMNIST, CIFAR-10 \citep{krizhevsky2009}, and the aloi dataset from the LIBSVM repository \citep{chang2011libsvm}. These three datasets respectively cover the cases of large number of data samples, large number of features, and large number of classes.They are summarized in Table~\ref{table_datasets}.

\begin{table}[h]
    \centering
\begin{tabular}{| c | c | c | c |}
    \hline
    Dataset  &  \# Features &  \# Classes &  \#  Examples \\ \hline
    aloi     & 128       & 1000     & 108,000 \\ \hline
    notMNIST & 784       & 10       & 500,000   \\ \hline
    CIFAR-10  & 3,072      & 10       & 60,000  \\ \hline
\end{tabular}
\caption{Summary of datasets}
\label{table_datasets}
\end{table}

For all experiments presented in this section, the underlying network topology is time-varying. In particular, the communication graph switches sequentially in a round robin manner within a pool of randomly generated doubly stochastic matrices. In addition, each dataset is divided into equally-sized parts, and distributed among the learner nodes in the network.

In the first experiment, we consider multi-class classification with logistic regression over a large distributed network consisting of $1500$ learner nodes. In this task, the learners observe a sequence of labeled examples $(\omega, z)$ taken from the aloi dataset, where $\omega \in \mathbb{R}^d$ denotes the feature vector, and $z \in \mathbb{R}$ represents the true label.  

For logistic regression, the loss associated with each data point is given by $f(x, (\omega_t, z_t)) = \log( 1 + \exp( -z_t x^T \omega_t )  )$. We use negative entropy as the regularization function $r(x)$. In addition, we set the batch size to $10$ data examples per online round, and use the step size of $\eta = 0.01$.

\begin{figure*}[t]
   
    \centering
   \includegraphics[scale =0.4]{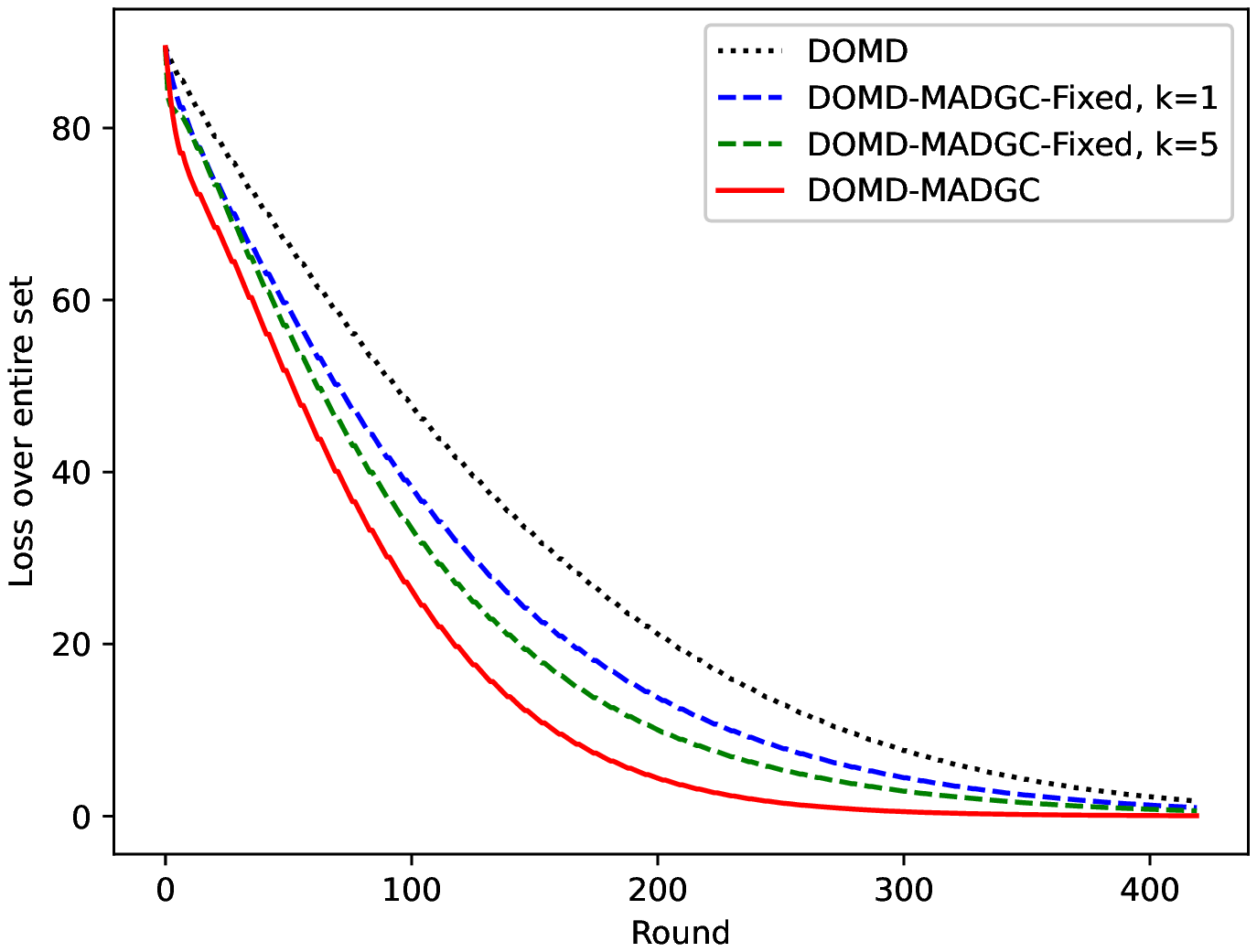}\caption{Logistic regression on aloi dataset.\label{f_loss_aloi}}

   \centering
   \hspace{0.1cm}

   \includegraphics[scale =0.4]{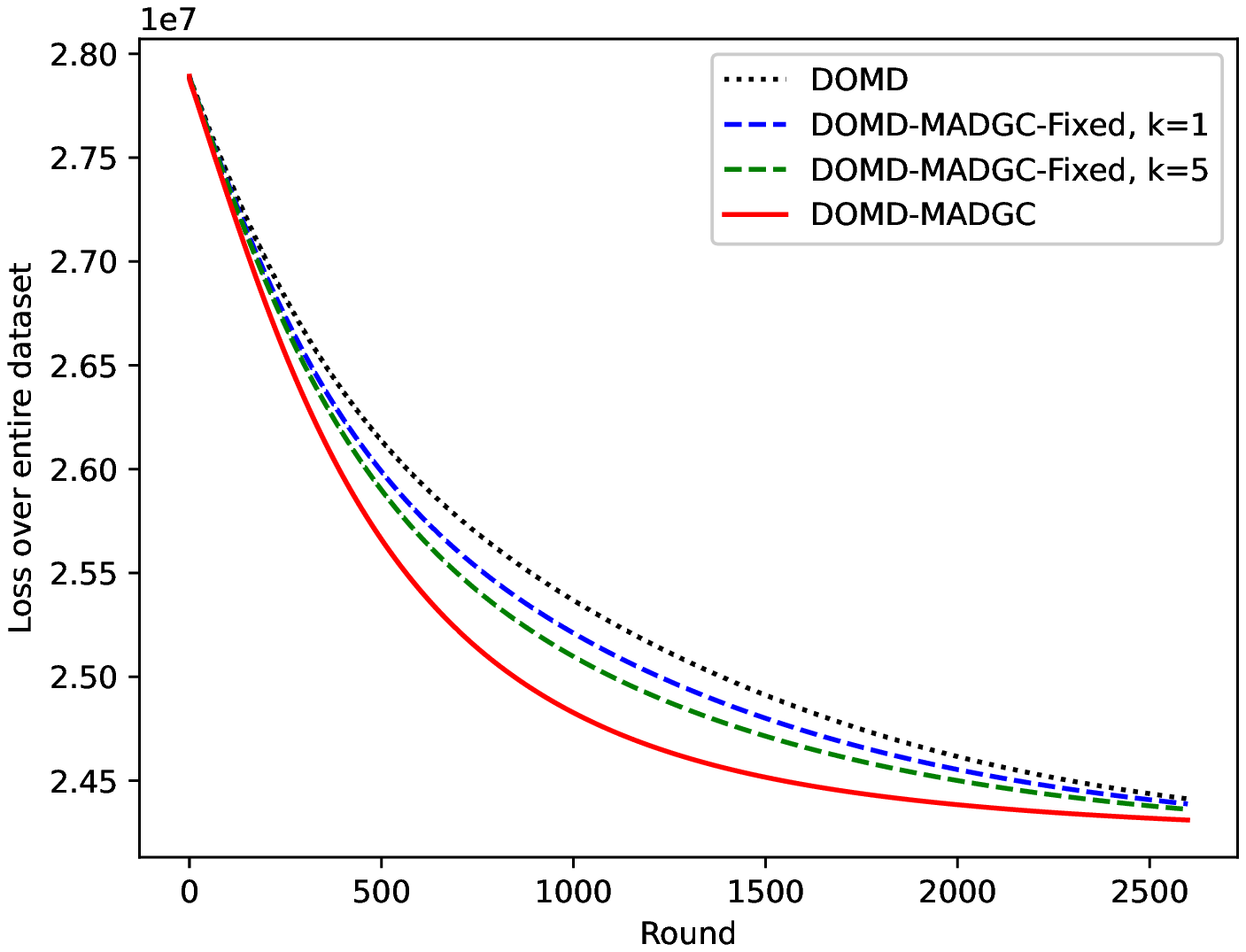}\caption{Ridge regression on notMNIST dataset.\label{f_loss_notMNIST}}
   
   \centering
   \hspace{0.1cm}

   \includegraphics[scale =0.4]{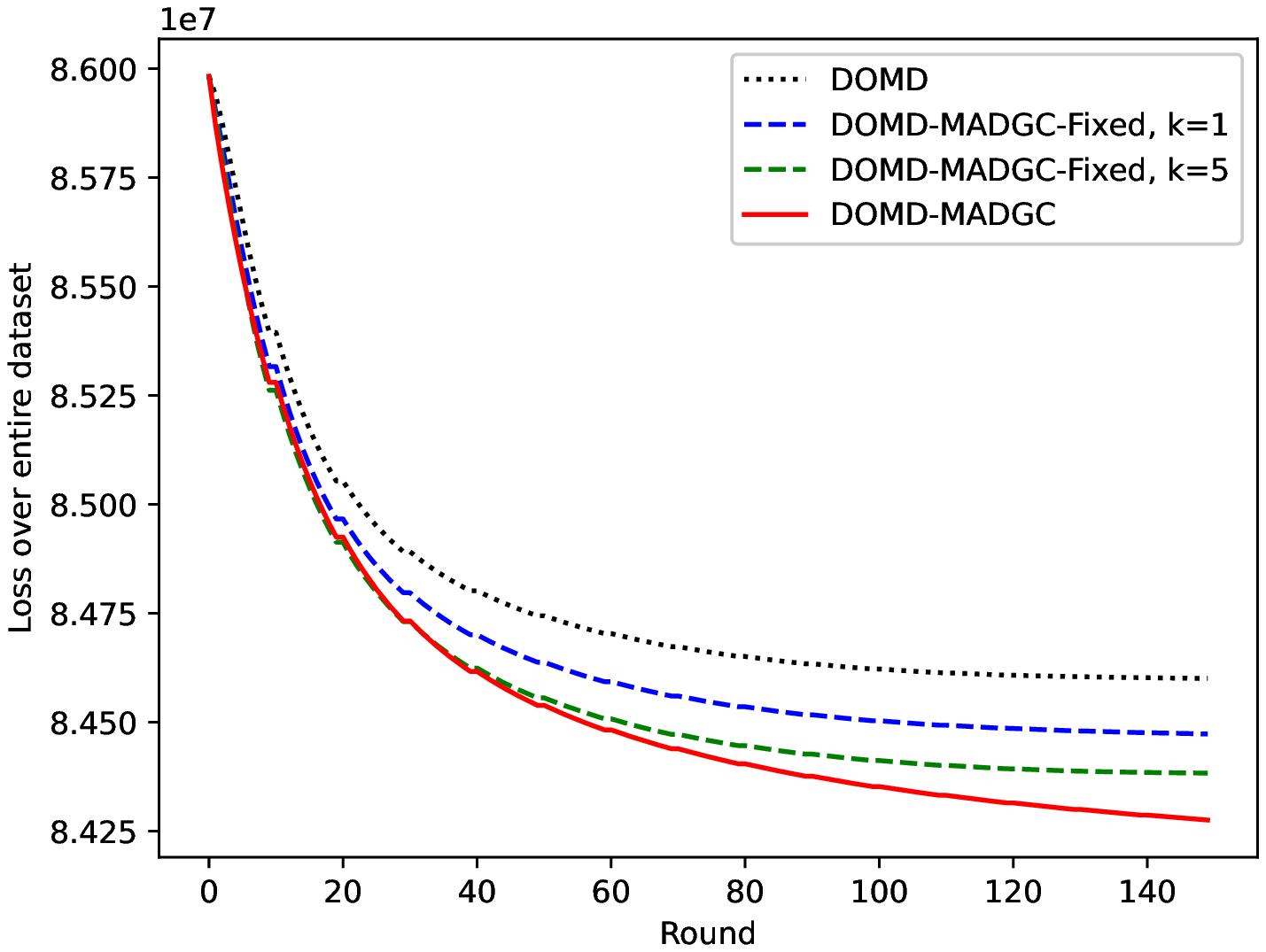}\caption{Ridge regression on CIFAR-10 dataset.\label{f_loss_cifar10}}

\end{figure*}

In Fig.~\ref{f_loss_aloi}, we compare the performance of DOMD-MADGC with fixed and increasing number of consensus rounds versus DOMD \citep{shahrampour2017distributed}.  Fig.~\ref{f_loss_aloi} shows that DOMD-MADGC, even with fixed and small $K = 1$ and $5$ outperforms the existing DOMD algorithm, which verifies the effectiveness of the proposed approach based on the communication of both local gradients and local decisions.  We also find that the overall loss of DOMD-MADGC with logarithmically increasing number of consensus iterations decreases faster than all alternatives.

Next, we consider the problem of ridge regression on the well-known notMNIST dataset, and over a network with $200$ distributed learners.  For this task, the loss associated with each data point is given by $f(x, (\omega_t, z_t)) = (x^T \omega_t - z_t)^2$. In the experiment, we set the batch size to $40$, set $\eta = 7 \times 10^{-4}$, and use the negative entropy as $r(x)$. Fig.~\ref{f_loss_notMNIST} shows the overall loss versus the number of online rounds. We find that DOMD-MADGC reduces the overall loss faster than the alternatives.  

We have also experimented on the CIFAR-10 dataset to solve the ridge regression problem in a distributed online manner. We consider a network of $500$ learner nodes, and set the batch size to $10$, and $\eta = 0.5$. Fig.~\ref{f_loss_cifar10} shows that DOMD-MADGC outperforms standard DOMD method in terms of the overall loss.

\end{document}